\documentclass[graybox]{svmult}

\RequirePackage{type1cm}        \RequirePackage{makeidx}         \RequirePackage{graphicx}        \RequirePackage{multicol}        \RequirePackage[bottom]{footmisc}\RequirePackage{amssymb}
\RequirePackage{url}

\DeclareFontFamily{U}{mathb}{\hyphenchar\font45}
\DeclareFontShape{U}{mathb}{m}{n}{ <-6> matha5 <6-7> matha6 <7-8>
mathb7 <8-9> mathb8 <9-10> mathb9 <10-12> mathb10 <12-> mathb12 }{}
\DeclareSymbolFont{mathb}{U}{mathb}{m}{n}

\DeclareMathAccent{\abxring}{0}{mathb}{"38}

\DeclareFontFamily{U}{mathb}{\hyphenchar\font45}
\DeclareFontShape{U}{mathb}{m}{n}{ <-6> matha5 <6-7> matha6 <7-8>
mathb7 <8-9> mathb8 <9-10> mathb9 <10-12> mathb10 <12-> mathb12 }{}
\DeclareSymbolFont{mathb}{U}{mathb}{m}{n}

\RequirePackage{newtxtext}       \usepackage{enumitem}

\RequirePackage{calrsfs}

\RequirePackage{amsfonts}

\RequirePackage[most]{tcolorbox}

\RequirePackage{bbm}

\RequirePackage{graphicx}

\RequirePackage{booktabs}

\RequirePackage{url}

\RequirePackage{hyperref}

\RequirePackage{color}

\RequirePackage{wasysym}
\RequirePackage{listings}
\RequirePackage{caption}
\RequirePackage{subcaption}

\RequirePackage{stmaryrd}

\RequirePackage{lineno}

\RequirePackage{xparse}

\RequirePackage{float}

\RequirePackage{multirow}

\RequirePackage{cellspace}

\RequirePackage{shellesc}

\RequirePackage[ruled,vlined,linesnumbered]{algorithm2e}

\SetCommentSty{xCommentSty}

\RequirePackage{tikz}
\usetikzlibrary{calc, automata, positioning, arrows, external, matrix, shapes.misc}
\tikzset{cross/.style={cross out,thick,draw=black,minimum size=2*(#1-\pgflinewidth), inner sep=0pt, outer sep=0pt},
cross/.default={2.5pt}}

\RequirePackage{pgfplots}
\pgfplotsset{compat=1.16}
\usepgfplotslibrary{fillbetween}
\usepgfplotslibrary{groupplots}
\usepgfplotslibrary{external}

\DeclareMathAlphabet{\pazocal}{OMS}{zplm}{m}{n}
\DeclareMathAlphabet\bpazocal{OMS}{cmsy}{b}{n}

\providecommand{\bR}{\mathbb{R}}

\providecommand{\vertiii}[1]{{\left\vert\kern-0.15ex\left\vert\kern-0.15ex\left\vert #1
    \right\vert\kern-0.15ex\right\vert\kern-0.15ex\right\vert}}

\NewDocumentCommand{\curlii}{sO{}m}
{
	\IfBooleanTF{#1}
    {\dgalext{#3}}
    {\dgalx[#2]{#3}}
}

\NewDocumentCommand{\dgalext}{m}{  \sbox0{    \mathsurround=0pt     $\left\{\vphantom{#1}\right.\kern-\nulldelimiterspace$  }  \sbox2{\{}  \ifdim\ht0=\ht2
    \{\kern-.625\wd2 \{#1\}\kern-.625\wd2 \}  \else
    \left\{\kern-.7\wd0\left\{#1\right\}\kern-.7\wd0\right\}  \fi
}

\NewDocumentCommand{\dgalx}{om}{  \sbox0{\mathsurround=0pt$#1\{$}  \sbox2{\{}  \ifdim\ht0=\ht2
    \{\kern-.625\wd2 \{#2\}\kern-.625\wd2 \}  \else
    \mathopen{#1\{\kern-.7\wd0 #1\{}
    #2
    \mathclose{#1\}\kern-.7\wd0 #1\}}
  \fi
}

\tikzset{
  partial ellipse/.style args={#1:#2:#3}{
    insert path={+ (#1:#3) arc (#1:#2:#3)}
  }
}

\providecommand{\pF}{\pazocal{F}}

\providecommand{\pV}{\pazocal{V}}

\providecommand{\to}{\widetilde{o}}

\definecolor{blackmy}{RGB}{38, 70, 83}

\definecolor{bluemy}{RGB}{39, 125, 161}
\definecolor{greenmy}{RGB}{42, 167, 143}
\definecolor{yellowmy}{RGB}{233, 196, 106}
\definecolor{brownmy}{RGB}{244, 162, 97}
\definecolor{redmy}{RGB}{249, 65, 68}

 \usepackage{xfrac}
\usepackage{color, colortbl}

\usepackage{lineno}

\makeatletter
\newcommand{\thickhline}{\noalign {\ifnum 0=`}\fi \hrule height 1pt
    \futurelet \reserved@a \@xhline
}
\newcolumntype{"}{@{\hskip\tabcolsep\vrule width 1pt\hskip\tabcolsep}}
\def\blfootnote{\gdef\@thefnmark{}\@footnotetext}
\makeatother

\makeindex

\begin{document}

\title*{Nonlinear Schwarz preconditioning for Quasi-Newton methods}

\author{Hardik Kothari}
\institute{Hardik Kothari  \at Universit\`{a} della Svizzera italiana, Switzerland, \email{hardik.kothari@usi.ch}
}
\maketitle
\abstract*{We propose the nonlinear restricted additive Schwarz (RAS) preconditioning strategy to improve the convergence speed of limited memory quasi-Newton (QN) methods.
We consider both ``left-preconditioning" and ``right-preconditioning" strategies.
As the application of the nonlinear preconditioning changes the standard gradients and Hessians to their preconditioned counterparts, the standard secant pairs cannot be used to approximate the preconditioned Hessians.
We discuss how to construct the secant pairs in the preconditioned QN framework.
Finally, we demonstrate the robustness and efficiency of the preconditioned QN methods using numerical experiments.}
\abstract{We propose the nonlinear restricted additive Schwarz (RAS) preconditioning strategy to improve the convergence speed of limited memory quasi-Newton (QN) methods.
We consider both ``left-preconditioning" and ``right-preconditioning" strategies.
As the application of the nonlinear preconditioning changes the standard gradients and Hessians to their preconditioned counterparts, the standard secant pairs cannot be used to approximate the preconditioned Hessians.
We discuss how to construct the secant pairs in the preconditioned QN framework.
Finally, we demonstrate the robustness and efficiency of the preconditioned QN methods using numerical experiments.}

\section{Introduction}
\label{sec:intro}
In this work, we consider a nonlinear preconditioning strategy for Quasi-Newton (QN) methods.
QN methods are a class of root-finding methods, where the full Jacobian is replaced with its approximation.
In the context of this work, we consider secant methods, which take into account a variable number of secant equations at each nonlinear iteration.
These types of methods are mostly used if the Jacobian of the nonlinear system is expensive to evaluate, requires more storage, or is simply unavailable.
Such scenarios are often encountered while solving coupled multiphysics problems that require higher-order discretization, inverse problems, optimal control problems, training of deep neural networks, etc.

To this aim, we consider the following abstract nonlinear minimization problem:
\begin{equation}
  \text{Find } x^{\ast} \in \pV \text{ that minimizes }\quad \Psi(x),
  \label{eq:min}
\end{equation}
where $\Psi:\pV \to \bR$ denotes a bounded, twice continuously differentiable objective function.
The objective function $\Psi$ is obtained by a finite element (FE) discretization of a nonlinear optimization problem, and $\pV$ denotes some FE space.
To solve \eqref{eq:min}, we can consider the first-order optimality condition for the function $\Psi(x)$, and then a nonlinear iterative method can be employed to find the root of the nonlinear equation $F(x^\ast) = 0$, where $F:\pV\to \pV'$ is defined as $F(\cdot) \equiv \nabla\Psi(\cdot)$.
We also note that the Hessian of the objective function $\nabla^{2} \Psi$ is equivalent to the Jacobian $F'$.
Depending on the properties of the objective function $\Psi$, multiple approaches can be considered to solve \eqref{eq:min}, for example, Newton's method and its variants; nonlinear Krylov methods; secant methods~\cite{jorgenocedal2000-04-27}.

Among all these methods, Newton's method is one of the most popular methods to solve such problems due to its locally quadratic convergence property.
However, its convergence might suffer if the objective function is highly nonlinear with locally stiff or unbalanced nonlinearities and if the initial guess is far from the solution.
In recent years, some nonlinear preconditioning strategies have been developed to accelerate the convergence of Newton's method, e.g.: Additive Schwarz Preconditioned Inexact Newton (ASPIN) \cite{cai_nonlinearly_2002};  Nonlinear Elimination Preconditioned Inexact Newton (NEPIN) \cite{cai_inexact_2011}; Restricted Additive Schwarz Preconditioned Exact Newton(RASPEN) \cite{dolean_nonlinear_2016}.
Similarly, in the context of optimization methods, nonlinear preconditioning strategies have been considered to improve the convergence of a nonlinear Krylov method~\cite{de_sterck_nonlinearly_2016} and a quasi-Newton (QN) method~\cite{de_sterck_nonlinearly_2018}.
To the best of our knowledge, unlike the ASPIN, NEPIN, and RASPEN methods, the nonlinear domain decomposition-based preconditioners have not yet been considered for Krylov methods and QN methods.

In this work, we apply the nonlinear Schwarz preconditioning strategies to accelerate the convergence of the standard QN method.
We explore the ``left" and ``right" nonlinear preconditioning strategies and discuss the necessary modifications to the QN framework.
Finally, we examine the efficiency of the preconditioned QN methods by means of some numerical experiments

\section{Preconditioned Quasi-Newton methods}
In this section, we discuss QN methods, nonlinear restricted additive Schwarz (NRAS) methods, and how to nonlinearly precondition QN methods.

\vspace{0.3cm}
\noindent\textbf{Quasi-Newton Methods: }
Quasi-Newton (QN) methods are quite popular in the optimization community, especially when the Hessian of the underlying minimization problem is unavailable or very expensive to evaluate.
In QN methods, the evaluation of the Hessian is replaced by its low-rank approximation.
This low-rank approximation of the Hessian is carried out using a secant condition.
At each iteration, the approximation of the Hessian $B$ is constructed using the information between subsequent iterations.
The approximate Hessian, $B^{(k+1)}$, satisfies the secant equation
\begin{equation}
  B^{(k+1)} s^{(k)} = y^{(k)},
\end{equation}
where  $s^{(k)} = x^{(k+1)}- x^{(k)}$ and $y^{(k)} = F(x^{(k+1)})- F(x^{(k)})$.
As the secant equation is not sufficient to uniquely determine the matrix $B$, additional constraints have to be imposed on $B$, which gives rise to different variants of the QN methods.
In this work, we consider two types of multi-secant methods, namely the Broyden-Fletcher-Goldfarb-Shanno (BFGS) method, and the Andersen acceleration (AA) method.
As one of the motivations of this work is reducing the memory footprint of the algorithm, the limited-memory variant of the BFGS method (L-BFGS), and of the AA method, becomes a natural choice.
These methods utilize only the $m$ pairs of the vectors $\{s^{(i)},y^{(i)}\}_{i=k-m}^{k-1}$ from the $m$ most recent iterations to construct the approximate Hessian.
We note that the original AA method is not proposed in the context of the optimization but its interpretation as a QN method is established in~\cite{fang_two_2009, zhang_globally_2020}.
The approximate Hessians obtained by the L-BFGS method and the type-I AA method (AA-I) at an iterate $k+1$ can be written in a compact matrix format in the following manner:
\begin{equation}
  \begin{aligned}
     & \text{(L-BFGS)} & B^{(k+1)} & = B_{0} - \begin{bmatrix} B_{0} S_{k} & Y_{k}\end{bmatrix} \begin{bmatrix} S_k B_0 S_k & L_k \\ L_k^\top & -D_k \end{bmatrix}^{-1} \begin{bmatrix} S_{k}^\top B_0 \\ Y_{k}^\top\end{bmatrix}, & \quad \\  & \text{(AA-I)}   & B^{(k+1)} & = I + (Y_{k} - S_{k})  ( S_{k}^\top S_{k})^{-1} S_{k}^\top.                                                                                                                                   &
    \label{eq:AA1}
  \end{aligned}
\end{equation}
Here, $S_k := [s^{(k-m)},\ldots,s^{(k-1)}]$, $Y_k := [y^{(k-m)},\ldots,y^{(k-1)}]$, $L_k$ and $D_k$ denote the strictly lower triangular, and the diagonal part of matrix $S_k^\top Y_k$, $B_0$ denotes some initial Hessian approximation.
In order to find the search direction $p^{(k)}$, we need the inverse of the approximate Hessians, which is generally obtained using the Sherman–Morrison–Woodbury formula. To accelerate the convergence speed of these methods, we propose to precondition the QN methods with an NRAS method.

\vspace{0.3cm}
\noindent\textbf{Nonlinear Restricted Additive Schwarz Methods: }
We consider a decomposition of the domain $\Omega$ into $n$ non-overlapping domains $\{\Omega_i\}_{i=1}^n$ and overlapping domains as $\{\Omega_i^\delta\}_{i=1}^n$, such that $\Omega_i \subset\Omega_i^\delta$, here $\delta$ denotes the size of the overlap.
The FE spaces associated with the overlapping domains are defined as $\{\pV_i^\delta\}_{i=1}^n$, $\pV_i^{\delta}\subset \pV$.
On these overlapping subspaces, we define the restriction and prolongation operators as $R^\delta_i:\pV\to\pV^\delta_i$ and $P^\delta_i:\pV^\delta_i\to \pV$, respectively.
We note that for $\delta=0$, the overlapping decomposition degenerates to a non-overlapping decomposition, i.e., $\Omega_i = \Omega_i^{0}$.
The prolongation operator on the non-overlapping subspaces is termed as \emph{restricted} prolongation operator, ${P_i^{0} : \pV^0_i \to \pV}$.
The overlapping and the non-overlapping decomposition of the subspaces ensures that the partition of unity is satisfied, e.g., $\sum_{i=1}^n P_i^{0} R_i^{\delta} = I$.

Now, we can define a local nonlinear minimization problem restricted to each overlapping subspace as follows.
For a given initial guess $x^\delta_i = R_i^\delta x^{(k)}$:
\begin{equation}
  \text{Find } x_i^{\ast}\in \pV_i^{\delta} \text{ that minimizes }\quad \Psi^\delta_i(x^\delta_i).
  \label{eq:local_min}
\end{equation}
Here, $\Psi^\delta_i :\pV_i^\delta \to \bR$ is the restriction of the objective function $\Psi$ to the subspace $\pV^{\delta}_i$.
Once the minimization problem is approximately solved on each subdomain, the global iterate is updated in the following manner
\begin{equation}
  x^{(k+1)} = x^{(k)} + \alpha^{(k)} \sum\nolimits_{i=1}^n P_i^{0} (x_i^\ast - R_i^\delta x^{(k)}).
  \label{eq:RAS}
\end{equation}
We note that the problem \eqref{eq:local_min} is solved on the overlapping subdomains, but the correction is accepted only on the non-overlapping part.
Furthermore, to construct a two-level variant of the NRAS method, we define a coarse space $\pV_0 \subset \pV$ and the restriction and the prolongation operators $R_0: \pV' \to \pV'_0 $ and $P_0 : \pV_0 \to \pV$, where $P_0^\top = R_0$.
Also, we define a projection operator $\Pi_0: \pV \to \pV_0$ to transfer the primal variables to the coarse level.
The objective function on the coarse level is defined as $\Psi_0:\pV_0 \to \bR$, which denotes a  discretization of the function $\Psi$ on the space $\pV_0$.
The coarse space plays an important role in the NRAS method, as it allows global communication between the subdomains and ensures the scalability of the algorithm.
In this work, the coarse-level objective function is defined in the spirit of the full approximation scheme (FAS) or the MG-Opt method~\cite{nash_multigrid_2000}.
The coarse-level function is constructed by adding a first-order consistency term, which is also called a ``defect" in the context of FAS.
Thus, the optimization problem on the coarse level is defined as follows.
For an initial guess $x_0 = \Pi_0 x^{(k)}$:
\begin{equation}
  \text{Find } x_0^\ast \text{ that minimizes }\quad \hat{\Psi}_0  (x_0) := \Psi_0(x_0) + \langle \delta g_0  ,x_0 \rangle
  \label{eq:coarse_min}
\end{equation}
where $\delta g_0 = R_0 \nabla\Psi(x^{(k)}) - \nabla \Psi_0(\Pi_0 x^{(k)})$ denotes the first-order consistency term.
Additionally, we employ a multiplicative variant of the coarse-level update, where we first approximately solve the problem on the coarse level and bootstrap the initial guess on the subdomains using the approximate solution from the coarse level.
The update step for the two-level NRAS is given as follows:
\begin{equation}
  x^{(k+1)} = x^{(k)} + \hat{\alpha} T_0(x^{(k)}) + \alpha \sum\nolimits_{i=1}^n P^0_i (x_i^\ast - R_i^\delta ( x^{(k)} + \hat{\alpha} T_0(x^{(k)})))
  \label{eq:TL_RAS}
\end{equation}
where $T_0(x^{(k)}) = P_0 (x_0^\ast - \Pi_0 x^{(k)})$ denotes the coarse-level correction.
We note that in \eqref{eq:TL_RAS}, $\hat{\alpha}$ and  $\alpha$ are computed using a line-search method, while $x_i^\ast$ and $x_0^\ast$ denote the approximate solutions of problems \eqref{eq:local_min} and \eqref{eq:coarse_min}, respectively.

\vspace{0.3cm}
\noindent\textbf{Nonlinear Preconditioning:}
In this section, we discuss strategies to nonlinearly precondition quasi-Newton methods.
Recall, we seek $x^\ast \in \pV$ such that $F(x^\ast) = 0$.
A nonlinear preconditioner $G$ of the residual function $F$ is defined such that the preconditioner $G$ approximates the inverse of the residual i.e., $G \approx F^{-1}$.
Practically, it is not possible to obtain such a preconditioning operator $G$ explicitly but, generally, such an operator can be defined implicitly as a fixed-point nonlinear iterative scheme, given as
\(
  x = G (x).
\)
The operator, $G$, can be applied to the nonlinear residual as either a ``left" or a ``right" preconditioner, which gives rise to two different nonlinearly preconditioned residuals
\begin{equation}
  \pF_L(x) = G_L(F(x)) = x - G(x), \quad \quad \quad \pF_R(x) = F(G_R(x))= F(G(x)).
\end{equation}
We remark that the left preconditioning operator is not equivalent to a fixed-point nonlinear iterative method $G_L \neq G$, while the right preconditioning operator is a fixed-point iteration scheme $G_R=G$.
The ASPIN and RASPEN methods are the ``left" preconditioned methods, where the nonlinear residual is first computed using a fixed-point method, and Newton's method is used to solve the equation $\pF_L(x)=0$.
The NEPIN method \cite{cai_inexact_2011}, nonlinear FETI-DP and BDDC methods \cite{klawonn_nonlinear_2017} are considered to be the ``right" preconditioned methods.

We define generic iterations for both types of preconditioning strategies.
The iteration for the preconditioned QN method can be achieved by replacing the residual $F$ with the preconditioned residual given as $\pF_{L/R}$.
For a given initial iterate $x^{(k)}$, we first compute $x^{(+)}$ using a NRAS method, i.e., ${x^{(+)} = G(x^{(k)})}$.
Once the preconditioning step has been carried out, we can define the iteration for the ``left-preconditioned" QN method as,
\begin{equation}
  x^{(k+1)} = x^{(k)} - \alpha^{(k)} \big(B_L^{(k)}\big)^{-1} \pF_L(x^{(k)}),  \ \ \text{where } \pF_L(x^{(k)}) = x^{(k)} - G(x^{(k)}).
  \label{eq:left}
\end{equation}
The update process for the ``right-preconditioned" QN method differs from the ``left-preconditioning" approach.
The iteration for the ``right-preconditioned" QN method is given as
\begin{equation}
  x^{(k+1)} = x^{(+)} - \alpha^{(k)} \big(B_R^{(k)}\big)^{-1} \pF_R(x^{(k)}),  \ \ \text{where } \pF_R(x^{(k)}) = F(G(x^{(k)})).
  \label{eq:right}
\end{equation}
In \eqref{eq:left} and \eqref{eq:right}, we compute $\alpha^{(k)}$ using a line-search method.
Here, $B_L^{(k)}$ and $B_R^{(k)}$ denote the approximation of the ``left" and ``right" preconditioned Hessians, respectively.
The QN method aims to approximate the Hessian of the underlying optimization function utilizing a set of vectors $\{s^k,y^k\}$. As we have preconditioned the QN method, we also have to change the underlying secant equation and corresponding secant pairs.
The corresponding secant equations for the ``left" and the ``right" preconditioned systems are now given as
\begin{equation}
  B_L^{(k+1)}  s_L^{(k)} = y_L^{(k)}, \qquad \qquad  B_R^{(k+1)}  s_R^{(k)} = y_R^{(k)}.
  \label{eq:prec_secant}
\end{equation}
From \eqref{eq:left} and \eqref{eq:right}, it is clear that $s_{L/R}^{(k)}$ at each  iteration are defined as corrections, which are given as
\begin{equation}
  s_L^{(k)} = x^{(k+1)} - x^{(k)}, \qquad \qquad  s_R^{(k)} = x^{(k+1)} - x^{(+)}.
  \label{eq:s_curr}
\end{equation}
Now, we focus our attention on the computation of $y_{L/R}^{(k)}$, which are defined as the difference between the preconditioned residuals
\begin{equation}
  y_L^{(k)} = \pF_L(x^{(k+1)}) - \pF_L(x^{(k)}), \qquad \qquad  y_R^{(k)} = F(x^{(k+1)}) - F(x^{(+)}).
  \label{eq:y_curr}
\end{equation}
We note that for the ``right" preconditioning approach,  the nonlinear preconditioner can be simplified as $\pF_R(x^{(k)}) = F(G(x^{(k)})) = F(x^{(+)})$, and the iteration in \eqref{eq:right} can be further simplified as
\[
  x^{(k+1)} = x^{(+)} - \alpha^{(k)} \big(B_R^{(k)}\big)^{-1} F(x^{(+)}).
\]
This update process can be interpreted as a half iteration, while the first half of the iteration is the preconditioning step $x^{(+)} = G(x^{(k)})$.
Hence, the ``right-preconditioned" QN method should only construct the approximation of the Hessian for the second half of the iteration.

A sketch of the nonlinearly preconditioned quasi-Newton method is provided in Algorithm \ref{alg:NRASQN}.

\begin{algorithm}[t]
  \caption{Nonlinearly Preconditioned QN method}
  \label{alg:NRASQN}
  \DontPrintSemicolon
  \SetKwInOut{Input}{Input}
  \SetKwInOut{Output}{Output}
  \KwData{$F:\pV \to \pV'$, $x^{(0)} \in \pV$, $k \mapsfrom 0$}
  \KwResult{$x^{(k)}$}
  \While {$\| F(x^{(k)}) \| \geqslant \epsilon_{\text{rtol}} \|F(x^{(0)})\|  $}
  {
    For given $x^{(k)}$, compute the preconditioned residual $\pF_{L/R} (x^{(k)})$ \\Compute direction using L-BFGS/AA-I approximation of preconditioned Hessian  $p_{L/R}^{(k)} \mapsfrom - \big(B_{L/R}^{{(k)}}\big)^{-1} \pF_{L/R}{x^{(k)}} $\\ Find $\alpha^{(k)}$ using a line-search algorithm \\ Update the iterate: $x^{(k+1)} \mapsfrom x^{(k)} + \alpha^{(k)} p_L^{(k)}$ or $x^{(k+1)} \mapsfrom x^{(+)} + \alpha^{(k)} p_R^{(k)}$ \\
    Compute $s_{L/R}^{(k)} $ as in \eqref{eq:s_curr} and $y_{L/R}^{(k)} $ using \eqref{eq:y_curr} \\
    Update the history of secant pairs $\{s_{L/R}^{(k)},{y}_{L/R}^{(k)}\}$\\
   Update $k \mapsfrom k+1$  \\ }
\end{algorithm}

\section{Numerical Experiments}
We investigate the performance of the nonlinearly preconditioned QN method through some numerical experiments.
To this aim, we consider a domain $\Omega = (0,1)^2$ with the boundary $\Gamma$.
The boundary $\Gamma$ is decomposed into four parts: top ($\Gamma_t = [0,1]\times \{1\}$), bottom ($\Gamma_b=[0,1]\times \{0\}$), left ($\Gamma_l= \{0\}\times [0,1]$) and right ($\Gamma_r = \{1\}\times[0,1]$).
We use the discretize-then-optimize approach, where the discretization is done with the first-order FE method using a mesh with $200 \times 200$ quadrilateral elements.
The coarse level is also constructed with the same approach, where a mesh with $10 \times 10$ elements is employed for discretization.

\vspace{0.3cm}
\noindent
\textbf{Minimal Surface:}
This experiment aims to find the surface of the minimal area described by the function $u$ by solving the following minimization problem:
\begin{equation}
  \begin{aligned}
    \min_{u \in H^1(\Omega)}  \Psi_M(u) = \int_\Omega \sqrt{(1 + \|\nabla u \|^2)} \ dx, \qquad \qquad \qquad \\
    \text{subject to } \begin{cases}
                         u = -0.5 \sin(2\pi x_2) \ \text{on } \Gamma_l, \ \ u = 0.5 \sin(2\pi x_2) \ \text{on } \Gamma_r, \\
                         u = -0.5 \sin(2\pi x_1) \ \text{on } \Gamma_b, \ \ u = 0.5 \sin(2\pi x_1) \ \text{on } \Gamma_t.
                       \end{cases}
  \end{aligned}
  \label{eq:min_surf}
\end{equation}

\vspace{0.3cm}
\noindent
\textbf{Setup for the solution methods:}
As we aim to study the behavior of the preconditioned QN method, we use a fixed configuration of the NRAS method.
The overlap for all experiments is prescribed as $\delta = 2$, and the domain $\Omega$ is decomposed into $8$ subdomains.
The partitioning of the mesh is carried out using the METIS library.
The preconditioned QN is terminated if one of these conditions is satisfied: $\|F(x^{(k)})\| < 10^{-7}$ or $\|F(x^{(k)})\|  < 10^{-6}\|F(x^{(0)})\|$.
The subdomain solvers in the NRAS method employ Newton's method, which terminates if  $\|F_i(x_i^{(k)})\| < 10^{-10}$ or $\|F_i(x_i^{(k)})\|  < 10^{-1}\|F_i(x_i^{(0)})\|$ is satisfied.
On the coarse level, we also employ Newton's method, which terminates if $\|F_0(x_0^{(k)})\| < 10^{-12}$ or $\|F_0(x_0^{(k)})\|  < 10^{-10}\|F_0(x_0^{(0)})\|$ is satisfied, also the maximum number of iterations is set to $5$.
The experiments are carried out using MATLAB on a system with an {Intel Core i9-9880H processor}, and $16$\,GB of memory.

\begin{table}[t]
  \centering
  \setlength{\tabcolsep}{2.5pt}
  \begin{tabular}{r"c|c"c|c"c|c"c|c"c|c}
    Memory     & \multicolumn{2}{c"}{m = 1} & \multicolumn{2}{c"}{m = 3} & \multicolumn{2}{c"}{m = 5} & \multicolumn{2}{c"}{m = 7} & \multicolumn{2}{c}{m = 10}                                                         \\ \thickhline
               & Time (s)                   & \# Iter                    & Time (s)                   & \# Iter                    & Time (s)                   & \# Iter     & Time (s) & \# Iter & Time (s) & \# Iter \\ \hline
    L-BFGS     & 698.16                     & 643                        & 720.01                     & 642                        & 699.53                     & 646         & 702.62   & 679     & 536.40   & 513     \\ \hline
    L-BFGS (L) & \emph{301.37}              & \emph{\ 25}                & 288.69                     & \ 23                       & 296.82                     & \ 24        & 288.61   & \ 23    & 300.00   & \ 25    \\
    L-BFGS (R) & 426.21                     & \ 36                       & 296.99                     & \ 22                       & \emph{278.68}              & \emph{\ 20} & 273.99   & \emph{\ 19}    & \emph{272.45}   & \emph{\ 20}    \\ \hline
AA-I (L)   & 350.60                     & \ 30                       & 285.56                     & \ 22                       & 284.20                     & \ 22        & 287.44   & \ 23    & 284.43   & \ 22    \\
    AA-I (R)   & 374.47                     & \ 36                       & \emph{274.40}              & \emph{\ 22}                & 281.74                     & \ 23        & \emph{269.98}   & \ 21    & 281.21   & \ 22    \\
  \end{tabular}
  \caption{Total number of iterations and the time to solution for different variants of the QN methods with and without preconditioning.}
  \label{tab:convergence_time_iter}
\end{table}
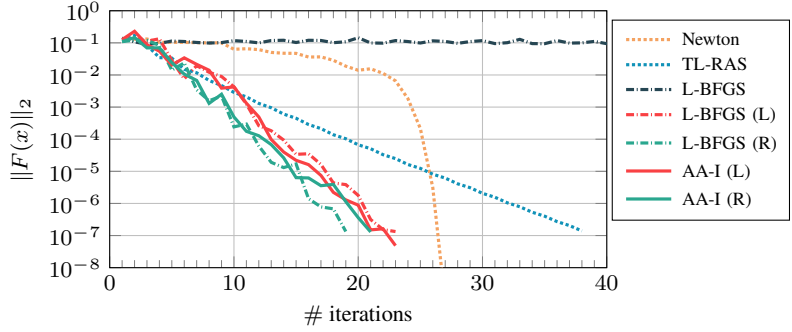
\begin{figure}[t]
  \centering
  \begin{tikzpicture}[]
    \begin{groupplot}[
        group style={
group size = 1 by 1,
x descriptions at=edge bottom,
horizontal sep=7pt,
          },
        width=0.7\textwidth,
        height=5cm,
        ymode=log,xlabel={$\#$ iterations}, minor tick num=9, grid=major,
        ymin=1e-8, ymax=1e0,
        xmax=40,
        xmin=0,label style={font=\small}, tick label style={font=\small}, legend style={font=\small}
      ]
\nextgroupplot[align=left, title={},ylabel={$\|F(x) \|_2$},ytick={1e0,1e-1,1e-2,1e-3,1e-4,1e-5,1e-6,1e-7,1e-8,1e-9}]
      \addplot[brownmy, densely dotted, very thick] table[x=iter, y=r_norm, col sep=comma] {./results/minsurf/Newton/Newton.csv};
      \label{pgfplots:Newton}
      \addplot[cyan!70!black, densely dotted, very thick] table[x=iter, y=r_norm, col sep=comma]{./results/minsurf/TL-RAS/SD8_OV2.csv};
      \label{pgfplots:RAS}
      \addplot[blackmy, densely dashdotted, very thick] table[x=iter, y=r_norm, col sep=comma] {./results/minsurf/LBFGS/Mem7.csv};
      \label{pgfplots:LBFGS}
      \addplot[redmy, densely dashdotted, very thick] table[x=iter, y=r_norm, col sep=comma]{./results/minsurf/LBFGS-L/SD8_OV2_Mem7.csv};
      \label{pgfplots:LBFGS-L}
      \addplot[greenmy, densely dashdotted, very thick] table[x=iter, y=r_norm, col sep=comma]{./results/minsurf/LBFGS-R/SD8_OV2_Mem7.csv};
      \label{pgfplots:LBFGS-R}
\addplot[redmy, very thick] table[x=iter, y=r_norm, col sep=comma]{./results/minsurf/AA-I-L/SD8_OV2_Mem7.csv};
      \label{pgfplots:AA1-L}
      \addplot[greenmy, very thick] table[x=iter, y=r_norm, col sep=comma]{./results/minsurf/AA-I-R/SD8_OV2_Mem7.csv};
      \label{pgfplots:AA1-R}
    \end{groupplot}
    \matrix [ fill=white,draw, matrix of nodes, anchor = north, node font=\scriptsize,
      column 1/.style={nodes={align=right,text width=0.6cm}},
      column 2/.style={nodes={align=left,text width=1.264cm}},
    ] at ($(group c1r1.south east) + (1.25, 3.3)$)
    {
      \ref{pgfplots:Newton}  & Newton     \\
      \ref{pgfplots:RAS}     & TL-RAS     \\
      \ref{pgfplots:LBFGS}   & L-BFGS     \\
      \ref{pgfplots:LBFGS-L} & L-BFGS (L) \\
      \ref{pgfplots:LBFGS-R} & L-BFGS (R) \\
\ref{pgfplots:AA1-L}   & AA-I (L)   \\
      \ref{pgfplots:AA1-R}   & AA-I (R)   \\
    };
  \end{tikzpicture}
  \caption{Convergence history of the standard QN methods, preconditioned QN methods, Newton's method, and the TL-NRAS method for solving the minimal surface problem. The QN methods are configured to use the last $7$ secant pairs.}
  \label{fig:convergence_history}
\end{figure}

\vspace{0.3cm}
\noindent
\textbf{Convergence study: }
In order to study the convergence behavior of the preconditioned QN method, the ``left" and the ``right" preconditioned variants of the L-BFGS methods and the AA-I method are considered.
For this numerical experiment, we decide to store $m$ pairs of secant vectors, where $m\in\{1,3,5,7,10\}$.
Table~\ref{tab:convergence_time_iter} depicts the time to solution and the required number of iterations to satisfy the termination criterion for different solution methods and different values of $m$.
We have included only preconditioned AA-I method and have excluded the vanilla AA-I method from our study.
\footnote{The AA-I method requires factorization of $S_k^\top Y_k$, which is not possible if the successive pairs of $\{s^{(k)},y^{(k)}\}$ are very similar.
To avoid such issues, one can construct the pairs in such a way that successive $s^{(k)}$ are orthogonal, but such modifications are out of the scope of this work.}
From Table~\ref{tab:convergence_time_iter}, it is clear that the preconditioned QN methods outperform the standard L-BFGS method both in terms of the number of iterations and the computational time.
Regardless of the number of stored secant pairs, the preconditioned L-BFGS methods and AA-I methods are two times faster than the L-BFGS method.
The preconditioned AA-I methods and the preconditioned L-BFGS methods have comparable performance. While the ``right" preconditioned L-BFGS methods outperform all other methods if more pairs of secant pairs are used.
Figure~\ref{fig:convergence_history} depicts the convergence history of the preconditioned QN methods, the two-level NRAS method, and Newton's method.
We can observe that the preconditioned QN method outperforms Newton's method and the L-BFGS method.
Also, we can see that the TL-NRAS method has linear convergence, and by employing a QN method as an outer solver we can reduce the number of required iterations in half.

 From the performed experiments, we can conclude that the proposed domain decomposition-based preconditioning strategy is quite robust both in the case of the L-BFGS method and the type-I AA method.
This works provides a promising future direction for problems when memory is a limiting factor, for example for solving the phase field fracture problems\cite{kopanivcakova2022nonlinear} or for the training of deep neural networks.

\vspace{0.1cm}
\footnotesize{
\noindent\textbf{Acknowledgements }
This work is supported by the Swiss National Science Foundation (SNSF) and the Deutsche Forschungsgemeinschaft for their through the project SPP 1962 ``Stress-Based Methods for Variational Inequalities in Solid Mechanics:~Finite Element Discretization and Solution by Hierarchical Optimization" [186407].
We also acknowledge the support of Platform for Advanced Scientific Computing through the project FraNetG:~Fracture Network Growth and SNSF through the project ML2 - Multilevel and Domain Decomposition Methods for Machine Learning [197041].
}
\bibliographystyle{spmpsci}
\bibliography{./extracted.bib}

\end{document}